\documentclass[11pt]{article}
\usepackage{amsfonts,latexsym,rawfonts,amsmath,amssymb,amsthm,graphicx}
\numberwithin{equation}{section}
\newtheorem{theorem}{Theorem}[section]

\newtheorem{thm}[theorem]{Theorem}

\newtheorem{defi}[theorem]{Definition}
\newtheorem{rem}[theorem]{Remark}
\def\s{\,\,\,\,}

\def\endproof{$\hfill\Box$\\}
\def\R{\mathbb{R}}
\def\H{\mathcal{H}}

\title{On biharmonic submanifolds in non-positively curved manifolds}
\author{Yong Luo}
\date{}
\begin{document}
\maketitle
\begin{abstract}
In the biharmonic submanifolds theory there is a generalized Chen's conjecture which states that biharmonic submanifolds in a Riemannian manifold with non-positive sectional curvature must be minimal. This conjecture turned out false by a counter example of Y. L. Ou and L. Tang in \cite{Ou-Ta}. However it remains interesting to find out sufficient conditions which guarantee this conjecture to be true. In this note we prove that:

 1. Any complete biharmonic submanifold (resp. hypersurface)  $(M, g)$ in a Riemannian manifold $(N, h)$ with non-positive sectional curvature (resp. Ricci curvature) which satisfies an integral condition: for some $p\in (0, +\infty)$, $\int_{M}|\vec{H}|^{p}du_g<+\infty,$
where $\vec{H}$ is the mean curvature vector field of $M\hookrightarrow N$, must be minimal. This generalizes the recent results due to N. Nakauchi and H. Urakawa in \cite{Na-Ur1} and \cite{Na-Ur2}.

 2. Any complete biharmonic submanifold (resp. hypersurface) in a Reimannian manifold of at most polynomial volume growth whose sectional curvature (resp. Ricci curvature) is non-positive must be minimal.

3. Any complete biharmonic submanifold (resp. hypersurface) in a non-positively curved manifold whose sectional curvature (resp. Ricci curvature) is smaller that $-\epsilon$ for some $\epsilon>0$ which satisfies that $\int_{B_\rho(x_0)}|\vec{H}|^{p+2}d\mu_g(p\geq0)$ is of at most polynomial growth of $\rho$, must be minimal.

 We also consider $\varepsilon$-superbiharmonic submanifolds  defined recently in \cite{Wh}  by G. Wheeler and prove similar results for $\varepsilon$-superbiharmonic submanifolds, which generalize the result in \cite{Wh}.

\end{abstract}

\section{Introduction}
Let $f: (M^m, g)\to (\R^{m+t}, h)$ be an isometric immersion from  a Riemannian manifold $M$ of dimension $m$ into an $m+t$-dimensional Euclidean space, where $t\geq1$ and $\vec{H}$  is the mean curvature vector field of $M$ in $\R^{m+t}$. $M$ is said to have harmonic  mean curvature vector field if (cf. \cite{Ch1})
 \begin{eqnarray}
 \Delta\vec{H}=0,
 \end{eqnarray}
 where $\Delta$ is the rough Laplacian on $(M, g)$.

 We see that a submanifold with harmonic mean curvature vector field is a natural generalization of minimal submanifolds.

Another natural generalization of minimal submanifolds is as follows. The intrinsic bi-energy of $f$ is defined by
\begin{eqnarray}
E_2(f)=\frac{1}{2}\int_M|\tau(f)|^2d\mu_g,
\end{eqnarray}
where $\tau(f)$ is the tension field of $f$ and $d\mu_g$ is the volume element on $(M, g)$. The critical points of the functional $E_2$ satisfy the following E-L equation ( see (\ref{equ3}))
\begin{eqnarray}\label{equ1}
-\Delta \vec{H}=\sum_{i=1}^mR^N(e_i, \vec{H})e_i,
\end{eqnarray}
where $R^N$ is the Riemann curvature tensor of the ambient manifold $N$ and $\{e_i, i=1,...,m\}$ is a local orthonormal frame of $M$.

A submanifold satisfying equation (\ref{equ1}) is called a biharmonic submanifold and we see that submanifolds with harmonic mean curvature vector fields in Euclidean spaces are biharmonic submanifolds. A natural question concerning biharmonic submanifolds is under what conditions they are minimal submanifolds. We have the following generalized Chen's conjecture (cf. \cite{Ca-Mo} \cite{Ca-Mo1} \cite{Ch} \cite{Ch1} \cite{Ou} \cite{Ou-Ta} and \cite{On} etc.).
\\\textbf{Generalized Chen's Conjecture.} Any biharmonic submanifold in a non-positively curved manifold is minimal.

Though it turned out that this conjecture is false by a counter example of Y. L. Ou and L. Tang (cf. \cite{Ou-Ta}), it remains interesting to find out sufficient conditions  which guarantee biharmonic submanifolds to be minimal. Recently  N. Nakauchi and H. Urakawa proved:
\begin{thm}[\cite{Na-Ur1} \cite{Na-Ur2}]
Assume that $(M, g)$ is a complete Riemannian manifold of dimension $m$ and $(N, h)$  a Riemannian manifold of dimension $m+t (t\geq1)$ whose sectional curvature (Ricci curvature if t=1) is non-positive. If $f: (M, g)\to (N, h)$ is a biharmonic isometric immersion with mean curvature vector field satisfying $$\int_M|\vec{H}|^2dM<+\infty,$$ then $f$ is minimal.
\end{thm}
In this paper, we generalize this result to the following:
\begin{thm}\label{thm}
Assume that $(M, g)$ is a complete Riemannian manifold of dimension $m$ and $(N, h)$  a Riemannian manifold of dimension $m +t (t\geq 1)$ whose sectional curvature is non-positive. If $f: (M, g)\to (N, h)$ is a biharmonic isometric immersion with mean curvature vector field satisfying  $$\int_M|\vec{H}|^pdM<+\infty,$$ for some $p\in (0, +\infty)$,  then $f$ is minimal. Furthermore, if $M$ is a hypersurface, then we only assume that $N$ has non-positive Ricci curvature.
\end{thm}
\begin{rem}
In a recent preprint \cite{Ma}, S. Maeta proved the above theorem for a complete biharmonic submanifold $M\hookrightarrow N$ with conditions:
\\1. $N$ has non-positive Ricci curvature, $M$ is a hypersurface and satisfies $\int_M|\vec{H}|^{\alpha}dM<+\infty$ for $1+\varepsilon<a<+\infty$, where $\varepsilon$ is a positive constant, or
\\2.  $N$ has non-positive sectional curvature and $M$ has general codimensions satisfying $\int_M|\vec{H}|^{\alpha}dM<+\infty$ for $1+\varepsilon<a<3-\varepsilon$, where $\varepsilon$ is a positive constant.

Clearly our assumption is weaker.
\end{rem}

We also find another sufficient condition to guarantee the Generalized Chen's conjecture to be true, to state our theorem, we first give a definition.
\begin{defi}\label{defi}
Let $(M, g)$ be a complete Riemannian manifold and $x_0\in M$. Then we say $(M, g)$ is of at most polynomial volume growth, if there exists a nonnegative integer $s$ such that
\begin{eqnarray}
Vol_g(B_\rho(x_0))\leq C\rho^s,
\end{eqnarray}
where $B_\rho(x_0)$ is the geodesic ball centered at $x_0$ with radius $\rho$ and  $C$ is a positive constat independent of $\rho$.
\end{defi}
\begin{rem}
We see that for any $\rho$ there is a constant $C_\rho$ which depends on $\rho$ such that $Vol_g(B_\rho(x_0))\leq C_\rho\rho^m$, where $m$ is the dimension of $M$. Hence the volume growth depends on the "blow up" rate of the metric at infinity.
\end{rem}
We have
\begin{thm}\label{thm0}
Assume that $(M, g)$ is a complete $m$-dimensional Riemannian manifold of at most polynomial volume growth and $(N, h)$  a Riemannian manifold of dimension $m + t (t\geq 1)$ whose sectional curvature is non-positive. Let $f:(M, g)\to (N, h)$ be an isometric immersion, then it is biharmonic if and only if it is minimal. Furthermore, if $M$ is a hypersurface, we only assume that $(N, h)$ has non-positive Ricci curvature.
\end{thm}
About the volume growth of a complete Riemannian manifold, the most celebrated result is the volume comparison theory due to Bishop and Gromov.  By Bishop-Gromov's volume comparison theorem, we see that manifolds with non-negative Ricci curvature have Euclidean volume growth (that is, $s=dim M$ in definition \ref{defi} ) and hence complete biharmonic submanifolds with nonnegative Ricci curvature are minimal, by theorem \ref{thm0}. But the assumption of nonnegative Ricci curvature seems  too strong. For volume comparison results of complete Riemannian manifold, we refer to a survey paper by G. Wei (see\cite{Wei}) and references therein.

We see all the above results give  restrictions on the submanifolds. In the following we weaken the assumption of integral of the mean curvature vector field by giving a condition on the target manifolds which guarantees the generalized Chen's conjecture to be true. In the following we let $B_\rho(x_0)$ be a geodesic ball of $M$ centered at $x_0$ of radius $\rho$.
\begin{thm}\label{thm4}
Assume that $(M,g)$ is a complete biharmonic submanifold  in a Riemannian manifold $(N,h)$ whose sectional curvature is smaller that $-\epsilon$ for some constant $\epsilon>0$  and $\int_{B_\rho(x_0)}|\vec{H}|^{p+2}d\mu_g(p\geq0)$ is of at most polynomial growth of $\rho$, then $M$ is a minimal submanifold. Furthermore, if $M$ is a hypersurface, we only assume that the Ricci curvature of $(N, h)$ is smaller than $-\epsilon$.
\end{thm}
We say a function $f: \R^+\to \R^+$ is of at most polynomial growth, if $f(\rho)\leq C(1+\rho^s)$ as $\rho\to \infty$, for some positive constant $C$ independent of $\rho$ and $s$ a positive integer.

If we look carefully at Ou and Tang's counter example (cf. \cite{Ou-Ta}) for the generalized
 Chen's conjecture we can see that their example is complete and the sectional curvature of the target manifold is negative and tends to zero at infinity.

\textbf{Question:} Is any complete biharmonic submanifold in a Riemannian manifold $(N,h)$ whose sectional curvature is smaller that $-\epsilon$ for some constant $\epsilon>0$ minimal?

For other type sufficient conditions which guarantee the generalized Chen's conjecture to be true, we refer the readers to papers  \cite{BMO} \cite{CMO} \cite{NUG}. In \cite{CMO}, Caddeo etc. proved that any biharmonic submanifold in hyperbolic 3-space $\H^3(-1)$ is minimal, and any
pseudo-umbilical biharmonic submanifold $M^m\subset \H^n(-1)$ with $m\neq4$ is minimal. It is
also shown in \cite{BMO} that any biharmonic hypersurface of $\H^n(-1)$ with at most two
distinct principal curvatures is minimal.

In a recent paper, G. Wheeler proposed a notion of $\varepsilon$-superbiharmonic submanifolds which is a generalization of submanifolds with harmonic mean curvature vector fields,  as follows:
\begin{defi}[\cite{Wh}]
Let $M$ be a submanifold in $N$ with metric $\langle.,.\rangle$, then we call $M$  an $\varepsilon$-superbiharmonic submanifold, if $$\langle\Delta \vec{H}, \vec{H}\rangle\geq (\varepsilon-1)|\nabla\vec{H}|^2,$$
where $\varepsilon\in [0,1]$ is a constant.
\end{defi}
\begin{rem}
When $\varepsilon=0$, we see that an $\varepsilon$-superbiharmonic submanifold satisfying $\Delta|\vec{H}|^2\geq0$, i.e. the square length of the mean curvature vector field of a 0-superbiharmonic submanifold is a superharmonic function.
\end{rem}
We have
\begin{thm}[\cite{Wh}]
Let $N^{m+t}$ be a complete Riemannian manifold. Suppose $f : M^m\to
N^{m+t}$ is a proper $\varepsilon$-superbiharmonic submanifold  for $\varepsilon>0$. Assume
in addition that $f$ satisfies the curvature growth condition
\begin{eqnarray}
\lim_{\rho\to\infty}\frac{1}{\rho^2}\int_{f^{-1}(B_\rho)}|\vec{H}|^2d\mu=0,
\end{eqnarray}
where $B_\rho$ is a geodesic ball of $N$ with radius $\rho$. Then $\vec{H}=0$ and $f$ is minimal. \end{thm}
We generalize this theorem to
\begin{thm}\label{thm2}
Let $N^{m+t}$ be a complete Riemannian manifold. Suppose $f : M^m\to
N^{m+t}$ is  a proper $\varepsilon$-superbiharmonic submanifold  for $\varepsilon>0$. Assume
in addition that $f$ satisfies the curvature growth condition
\begin{eqnarray}
\lim_{\rho\to\infty}\frac{1}{\rho^2}\int_{f^{-1}(B_\rho)}|\vec{H}|^{2+a}d\mu=0,
\end{eqnarray}
where $B_\rho$ is a geodesic ball of $N$ with radius $\rho$ and $a\geq0$ is a constant. Then $\vec{H}=0$ and $f$ is minimal.
\end{thm}
\begin{rem}
We say a map is proper, if the preimage of any compact subset of the ambient manifold is compact. Note that properness implies completeness of the metric of the domain manifold.
\end{rem}
If we drop the assumption of properness of the immersion,  we have
\begin{thm}\label{thm3}
Suppose $f : M^m\to
N^{m+t}$ is  a complete $\varepsilon$-superbiharmonic submanifold  in $N^{m+t}$ for $\varepsilon>0$. Assume
in addition that $f$ satisfies the condition
\begin{eqnarray}
\int_M|\vec{H}|^{2+a}d\mu<\infty,
\end{eqnarray}
where  $a\geq0$ is a constant. Then $\vec{H}=0$ and $f$ is minimal.
\end{thm}
Several words on the proof: The main argument in our proof is integral by parts. But compared with some other results in this line, instead of using integral by parts to the biharmonic submanifold equation, we integrate over the tangential part of the biharmonic submanifolds equation, and when we integrate by parts, we chose test functions more delicately  and  use the Young's inequality in a very subtle way.

\textbf{Organization.} In section 2 we give some preliminaries on submanifolds theory, harmonic and biharmonic maps. Theorems \ref{thm}, \ref{thm0}, \ref{thm4}, \ref{thm2}, \ref{thm3} are proved in section 3.

\section{Preliminaries}
\subsection{Submanifolds}
Assume that $f:M\to (N, h)$ is an immersion from a manifold of dimension $m$ to a Riemannian manifold $(N, h)$ of dimension $m+t$. Then $M$ inherits a
Riemannian metric from $(N, h)$ by $g(X, Y):=h(df(X) , df(Y)) $ for any $X, Y\in TM$
and a volume form by $d\mu_g:=\sqrt{\det g}dx$.
The second fundamental form of $M\hookrightarrow N$, $ B: TM\otimes TM\to NM$, is defined by
$$B(X, Y):=D_XY-\nabla_XY,$$
for any $X, Y\in TM$, where $D$ is the covariant derivative with respect to the Levi-Civita connection on $N$, $\nabla$ is the Levi-Civita connection on $M$ with respect to the induced metric and $NM$ is the normal bundle of $M$. For any normal vector field $\eta$ the Weingarten map associated with $\eta$, $A_\eta: TM\to TM$ is defined by
\begin{eqnarray}
D_X\eta=-A_\eta X+\nabla^\bot_X\eta,
\end{eqnarray}
where $\nabla^\bot$ is the normal connection and as is well known that $B$ and $A$ are related by
\begin{eqnarray}
\langle B(X, Y), \eta\rangle=\langle A_\eta X, Y\rangle.
\end{eqnarray}

For any $x\in M$, let $\{e_1, e_2,...,e_m, e_{m+1},...,e_{m+t}\}$ be a local orthonormal basis of $N$ such that $\{e_1,...,e_m\}$ is an orthonormal basis of $T_xM$. Then $B$ is decomposed  at $x$ as
$$B(X,Y)=\sum_{\alpha=m+1}^{m+t}B_\alpha(X, Y)e_\alpha.$$
The mean curvature vector field is defined as
\begin{eqnarray}
\vec{H}:=\frac{1}{m}\sum_{i=1}^mB(e_i, e_i)=\sum_{\alpha=m+1}^{m+t}H_\alpha e_\alpha,
\end{eqnarray}
where
\begin{eqnarray}
H_\alpha:=\frac{1}{m}\sum_{i=1}^mB_\alpha(e_i, e_i).
\end{eqnarray}
\subsection{Harmonic and biharmonic maps}
Let $f:(M, g)\to (N, h)$ be a map from a Riemannian manifold $(M, g)$ to a Riemannian manifold $(N, h)$. The Dirichlet energy of $f$ is defined by
$$E(f):=\frac{1}{2}\int_M|df|^2d\mu_g.$$
The E-L equation of $E$ is
$$\tau(f)=\rm tr\nabla df=0,$$
where $\tau(f)$ is called the tension field of $f$. A map satisfying this E-L equation is called a harmonic map.

To generalize the notion of harmonic maps, a natural way was proposed by J. Eells and L. Lemaire (\cite{EL}) in 1983. They considered the bi-energy $E_2(f):=\frac{1}{2}\int_M|\tau(f)|^2d\mu_g$. Critical points of the bi-energy are called biharmonic maps and we see that all harmonic maps are minimizers of the bi-energy functional $E_2$. In 1986, G. Y. Jiang (\cite{Ji}) calculated the first and second variational formulas of the bi-energy functional. The E-L equation of $E_2$ is
\begin{eqnarray}\label{equ2}
\tau_2(f)=-\Delta \tau(f)-\sum_{i=1}^mR^N(\tau(f), df(e_i))df(e_i)=0,
\end{eqnarray}
where $\Delta$ is the rough Laplacian on $(M, g)$, $\{e_i, i=1,...,m\}$ is a local orthonormal frame on $M$ and $R^N$ is the Riemann curvature tensor of $(N, h)$.

A map $f:(M, g)\to (N, h)$ satisfying equation (\ref{equ2}) is called a biharmonic map. Further if $f$ is an isometric immersion, it is called a biharmonic submanifold.

 Assume that $f$ is an isomeric immersion. we see $\{df(e_i)\}$ is a local orthonormal frame of $M$. In addition, for any $X, Y\in TM$,
 \begin{eqnarray*}
 \nabla d(X, Y)=\nabla^f_X(df(Y))-df(\nabla_XY),
 \end{eqnarray*}
 where $\nabla^f$ is the connection on the pullback bundle $f^{-1}TN$, whose fibre at a point $x\in M$ is $T_{f(x)}N=TM\oplus NM$. Here $TM$ is the tangent bundle and $NM$ is the normal bundle of $M$. Thus we see that $\nabla^f_X(df(Y))-df(\nabla_XY)=B(X, Y)$, and so
 \begin{eqnarray}
 \tau(f)=tr\nabla df=tr B=m\vec{H}.
 \end{eqnarray}
Therefore a biharmonic submanifold satisfying the following equation:
\begin{eqnarray}\label{equ3}
-\Delta \vec{H}-\sum_{i=1}^mR^N( e_i, \vec{H})e_i=0,
\end{eqnarray}
where $\Delta$ is the rough Laplacian on $(M, g)$, $\{e_i, i=1,...,m\}$ is a local orthonormal frame on $M$ and $R^N$ is the Riemann curvature tensor of $(N, h)$.

By decomposing $-\Delta \vec{H}-\sum_{i=1}^mR^N( e_i, \vec{H})e_i$ into its tangential and normal parts, we see that a submanifold is biharmonic if and only if it satisfies (cf. \cite{BMO1})
\begin{eqnarray}
\Delta^\bot\vec{H}-\sum_{i=1}^mB(A_{\vec{H}}e_i, e_i)+\sum_{i=1}^m(R^N( e_i, \vec{H})e_i)^\bot&=&0,\label{equ4}
\\m\nabla|\vec{H}|^2+4\sum_{i=1}^mA_{\nabla^\bot_{e_i} \vec{H}}e_i-\sum_{i=1}^m(R^N( e_i, \vec{H})e_i)^T&=&0,
\end{eqnarray}
where $(R^N( e_i, \vec{H})e_i)^\bot, (R^N( e_i, \vec{H})e_i)^T $ denote the normal and tangential parts of $R^N( e_i, \vec{H})e_i$ respectively. In particular, if $N$ is a space form of constant sectional curvature $c$, then a submanifold is biharmonic if and only if (cf. \cite{Ch2, On})
\begin{eqnarray}
\Delta^\bot\vec{H}-\sum_{i=1}^mB(A_{\vec{H}}e_i, e_i)+cm\vec{H}&=&0,
\\m\nabla|\vec{H}|^2+4\sum_{i=1}^mA_{\nabla^\bot_{e_i} \vec{H}}e_i&=&0.
\end{eqnarray}
\section{Proof of theorems}
In this section we prove our theorems.
\subsection{Proof of theorem \ref{thm}.}
\proof From equation $(\ref{equ4})$ we see that
\begin{eqnarray*}
\Delta|\vec{H}|^2&=&2|\nabla^\bot\vec{H}|^2+2\langle\vec{H}, \Delta^\bot\vec{H}\rangle
\\&=&2|\nabla^\bot\vec{H}|^2+2\sum_{i=1}^m\langle B(A_{\vec{H}}e_i, e_i), \vec{H}\rangle-2\sum_{i=1}^m\langle R^N( e_i, \vec{H})e_i, \vec{H}\rangle
\\&=&2|\nabla^\bot\vec{H}|^2+2\sum_{i=1}^m\langle A_{\vec{H}}e_i, A_{\vec{H}}e_i\rangle-2\sum_{i=1}^m\langle R^N( e_i, \vec{H})e_i, \vec{H}\rangle.
\end{eqnarray*}
Obviously if $N$ has non-positive sectional curvature, $-\sum_{i=1}^m\langle R^N( e_i, \vec{H})e_i, \vec{H}\rangle\geq0$. Furthermore, if $M$ is a hypersurface, we see that $-\sum_{i=1}^m\langle R^N( e_i, \vec{H})e_i, \vec{H}\rangle=-Ric^N(\vec{H}, \vec{H})\geq0.$
Therefore
\begin{eqnarray}\label{ineluo}
\Delta|\vec{H}|^2&\geq&2|\nabla^\bot\vec{H}|^2+2\sum_{i=1}^m\langle A_{\vec{H}}e_i, A_{\vec{H}}e_i\rangle \nonumber
\\&\geq&2|\nabla^\bot\vec{H}|^2+2m|\vec{H}|^4,
\end{eqnarray}
where in the last inequality we used $\sum_{i=1}^m\langle A_{\vec{H}}e_i, A_{\vec{H}}e_i\rangle\geq m|\vec{H}|^4$, which can be seem as follows: Let $x\in M$, when $\vec{H}(x)=0$, we are done. If $\vec{H}(x)\neq 0$, set $e_{m+t}=\frac{\vec{H}}{|\vec{H}|}$, then $\vec{H}(x)=H_{m+t}(x)e_{m+t}$ and $|\vec{H}|^2=H^2_{m+t}$. Hence we have at $x$,
\begin{eqnarray*}
\sum_{i=1}^m\langle A_{\vec{H}}e_i, A_{\vec{H}}e_i\rangle&=&H^2_{m+t}\sum_{i=1}^m\langle A_{e_{m+t}}e_i, A_{e_{m+t}}e_i\rangle
\\&=&H^2_{m+t}|B_{m+t}|_h^2
\\&\geq&mH^4_{m+t}
\\&=&m|\vec{H}|^4.
\end{eqnarray*}

Let $\gamma: M\to \R^+$ be a cut off function such that
$$\gamma=1\s on \s B_{\rho}, \gamma=0\s on\s M\setminus B_{2\rho},
\s and \s |\nabla\gamma|\leq\frac{C}{\rho},$$
for some constant $C$ independent of $\rho$. Here $B_\rho$ is a geodesic ball of radius $\rho$ on $M$.

By integral by parts and (\ref{ineluo}) one gets
\begin{eqnarray}\label{ine31}
&&-\int_M\nabla\gamma^b|\vec{H}|^a\nabla |\vec{H}|^2d\mu_g \nonumber
\\&=&\int_M\gamma^b|\vec{H}|^a\Delta|\vec{H}|^2d\mu_g  \nonumber
\\&\geq&2\int_M\gamma^b|\vec{H}|^a|\nabla^\bot\vec{H}|^2d\mu_g
+2m\int_M\gamma^b|\vec{H}|^{a+4}d\mu_g ,
\end{eqnarray}
where $a, b$ are positive constants to be determined later.

On the other hand one has
\begin{eqnarray}\label{ine41}
&&-\int_M\nabla\gamma^b|\vec{H}|^a\nabla |\vec{H}|^2d\mu_g \nonumber
\\&=&-2\int_Mb\gamma^{b-1}\nabla\gamma|\vec{H}|^a\langle\nabla^\bot\vec{H}, \vec{H}\rangle d\mu_g
-2a\int_M\gamma^b|\vec{H}|^{a-2}\langle\nabla^\bot\vec{H}, \vec{H}\rangle^2 d\mu_g \nonumber
\\&\leq&-2\int_Mb\gamma^{b-1}\nabla\gamma|\vec{H}|^a\langle\nabla^\bot\vec{H}, \vec{H}\rangle d\mu_g .
\end{eqnarray}
From (\ref{ine31})-(\ref{ine41}) we obtain,
\begin{eqnarray*}\label{ine5}
&&2\int_M\gamma^b|\vec{H}|^a|\nabla^\bot\vec{H}|^2d\mu_g
+2m\int_M\gamma^b|\vec{H}|^{a+4}d\mu_g
\\&\leq&-2\int_Mb\gamma^{b-1}\nabla\gamma|\vec{H}|^a\langle\nabla^\bot\vec{H}, \vec{H}\rangle d\mu_g .
\end{eqnarray*}
Now let $b=a+4$ we get
\begin{eqnarray}\label{ine7}
&&2\int_M\gamma^{a+4}|\vec{H}|^a|\nabla^\bot\vec{H}|^2d\mu_g  \nonumber
+2m\int_M\gamma^{a+4}|\vec{H}|^{a+4}d\mu_g
\\&\leq&-2\int_M(a+4)\gamma^{a+3}\nabla\gamma|\vec{H}|^a\langle\nabla^\bot\vec{H}, \vec{H}\rangle d\mu_g .\nonumber
\\&\leq&\int_M\gamma^{a+4}|\vec{H}|^a|\nabla^\bot\vec{H}|^2d\mu_g +
(a+4)^2\int_M\gamma^{a+2}|\vec{H}|^{a+2}|\nabla\gamma|^2d\mu_g  \nonumber
\\&=&\int_M\gamma^{a+4}|\vec{H}|^a|\nabla^\bot\vec{H}|^2d\mu_g +
(a+4)^2\int_M\gamma^q|\vec{H}|^q\gamma^{a+2-q}|\vec{H}|^{a+2-q}|\nabla\gamma|^2d\mu_g ,
\end{eqnarray}
where $q$ is a constant belonging to $(0, a+2)$.

By Young's inequality,
\begin{eqnarray}\label{ine8}
&&(a+4)^2\int_M\gamma^q|\vec{H}|^q\gamma^{a+2-q}|\vec{H}|^{a+2-q}|\nabla\gamma|^2d\mu_g \nonumber
\\&\leq&\int_M\gamma^{a+4}|\vec{H}|^{a+4}d\mu_g \nonumber
\\&+&C(a, q)\int_M\gamma^{(a+2-q)\frac{a+4}{a+4-q}}|\vec{H}|^{(a+2-q)\frac{a+4}{a+4-q}}
|\nabla\gamma|^\frac{2(a+4)}{a+4-q}d\mu_g \nonumber
\\&=&\int_M\gamma^{a+4}|\vec{H}|^{a+4}d\mu_g +C(a, q)\int_M\gamma^{a+4-\frac{2(a+4)}{a+4-q}}|\vec{H}|^{a+4-\frac{2(a+4)}{a+4-q}}
|\nabla\gamma|^\frac{2(a+4)}{a+4-q}d\mu_g
\end{eqnarray}
where $C(a, q)$ is a constant depending on $a, q$.

From (\ref{ine7})-(\ref{ine8}) we obtain
\begin{eqnarray*}
&&\int_M\gamma^{a+4}|\vec{H}|^a|\nabla^\bot\vec{H}|^2d\mu_g +(2m-1)\int_M\gamma^{a+4}|\vec{H}|^{a+4}d\mu_g
\\&\leq&C(a, q)\frac{1}{\rho^\frac{2(a+4)}{a+4-q}}\int_M\gamma^{a+4-\frac{2(a+4)}
{a+4-q}}|\vec{H}|^{a+4-\frac{2(a+4)}{a+4-q}}d\mu_g
\\&\leq&C(a, q)\frac{1}{\rho^\frac{2(a+4)}{a+4-q}}\int_M|\vec{H}|^{a+4-\frac{2(a+4)}{a+4-q}}d\mu_g
\end{eqnarray*}
Note that when $q$ varies from $0$ to $a+2$, $a+4-\frac{2(a+4)}{a+4-q}$ varies from $a+2$ to $0$.
Now by assumption we assume that for some constant $p\in(0, +\infty)$
$$\int_M|\vec{H}|^pd\mu_g <+\infty,$$
then we can choose $q$ and $a$ such that $a+4-\frac{2(a+4)}{a+4-q}=p,$ and so
\begin{eqnarray*}
&&\int_M\gamma^{a+4}|\vec{H}|^a|\nabla^\bot\vec{H}|^2d\mu_g +(2m-1)\int_M\gamma^{a+4}|\vec{H}|^{a+4}d\mu_g
\\&\leq&C(a, q)\frac{1}{\rho^\frac{2(a+4)}{a+4-q}}\int_M|\vec{H}|^pd\mu_g .
\end{eqnarray*}
Let $\rho\to \infty$ we have
$$\int_M\gamma^{a+4}|\vec{H}|^a|\nabla^\bot\vec{H}|^2d\mu_g
+(2m-1)\int_M\gamma^{a+4}|\vec{H}|^{a+4}d\mu_g =0,$$
which implies that $\vec{H}=0$.

This completes the proof of theorem \ref{thm}. \endproof

\subsection{Proof of theorem \ref{thm0}.}
The proof of theorem \ref{thm0} is quite similar to the proof of theorem \ref{thm}. But we still give all the details.
\proof From the beginning of the proof of theorem \ref{thm}, we see that if $N$ has non-positive sectional curvature or $M$ is a hypersurface and $N$ has non-positive Ricci curvature, then $-\sum_{i=1}^m\langle R^N( e_i, \vec{H})e_i, \vec{H}\rangle\geq0$.

From equation $(\ref{equ4})$ we see that
\begin{eqnarray}\label{ineluo1}
\Delta|\vec{H}|^2&=&2|\nabla^\bot\vec{H}|^2+2\langle\vec{H}, \Delta^\bot\vec{H}\rangle \nonumber
\\&=&2|\nabla^\bot\vec{H}|^2+2\sum_{i=1}^m\langle B(A_{\vec{H}}e_i, e_i), \nonumber \vec{H}\rangle-2\sum_{i=1}^m\langle R^N( e_i, \vec{H})e_i, \vec{H}\rangle
\\&\geq&2|\nabla^\bot\vec{H}|^2+2m|\vec{H}|^4.
\end{eqnarray}
As before,  let $\gamma: M\to \R^+$ be a cut off function such that
$$\gamma=1\s on \s B_{\rho}, \gamma=0\s on\s M\setminus B_{2\rho},
\s and \s |\nabla\gamma|\leq\frac{C}{\rho},$$
for some constant $C$ independent of $\rho$. Here $B_\rho$ is a geodesic ball of radius $\rho$ on $M$.

By integral by parts and (\ref{ineluo1}) one gets
\begin{eqnarray}\label{ine3}
&&-\int_M\nabla\gamma^b|\vec{H}|^a\nabla |\vec{H}|^2d\mu_g  \nonumber
\\&=&\int_M\gamma^b|\vec{H}|^a\Delta|\vec{H}|^2d\mu_g \nonumber
\\&\geq&2\int_M\gamma^b|\vec{H}|^a|\nabla^\bot\vec{H}|^2d\mu_g
+2m\int_M\gamma^b|\vec{H}|^{a+4}d\mu_g ,
\end{eqnarray}
where $a, b$ are positive constants to be determined later.

On the other hand one has
\begin{eqnarray}\label{ine4}
&&-\int_M\nabla\gamma^b|\vec{H}|^a\nabla |\vec{H}|^2d\mu_g  \nonumber
\\&=&-2\int_Mb\gamma^{b-1}\nabla\gamma|\vec{H}|^a\langle\nabla^\bot\vec{H}, \vec{H}\rangle d\mu_g
-2a\int_M\gamma^b|\vec{H}|^{a-2}\langle\nabla^\bot\vec{H}, \vec{H}\rangle^2 d\mu_g  \nonumber
\\&\leq&-2\int_Mb\gamma^{b-1}\nabla\gamma|\vec{H}|^a\langle\nabla^\bot\vec{H}, \vec{H}\rangle d\mu_g .
\end{eqnarray}
From (\ref{ine3})-(\ref{ine4}) we obtain,
\begin{eqnarray}\label{ine5}
2\int_M\gamma^b|\vec{H}|^a|\nabla^\bot\vec{H}|^2d\mu_g
+2m\int_M\gamma^b|\vec{H}|^{a+4}dM\leq
-2\int_Mb\gamma^{b-1}\nabla\gamma|\vec{H}|^a\langle\nabla^\bot\vec{H}, \vec{H}\rangle d\mu_g .
\end{eqnarray}
Now let $\frac{a+4}{a+2}=\frac{b}{b-2}$, i.e. $b=a+4$. Using Young's inequality, we have
\begin{eqnarray}\label{ineluo2}
&&-2\int_M(a+4)\gamma^{a+3}\nabla\gamma|\vec{H}|^a\langle\nabla^\bot\vec{H}, \vec{H}\rangle d\mu_g \nonumber
\\&\leq&\int_M\gamma^{a+4}|\vec{H}|^a|\nabla^\bot\vec{H}|^2d\mu_g
+\int_M\gamma^{a+4}|\vec{H}|^{a+4}d\mu_g
+C(a)\int_M|\nabla\gamma|^\frac{a+4}{2}d\mu_g,
\end{eqnarray}
where $C(a)$ is a constant depending on $a$.

From (\ref{ine5})-(\ref{ineluo2}) we obtain
\begin{eqnarray}
&&\int_M\gamma^{a+4}|\vec{H}|^a|\nabla^\bot\vec{H}|^2d\mu_g  \nonumber
+(2m-1)\int_M\gamma^{a+4}|\vec{H}|^{a+4}d\mu_g
\\&\leq&C(a)\int_M|\nabla\gamma|^\frac{a+4}{2}d\mu_g  \nonumber
\\&\leq&C(a)\frac{1}{\rho^{\frac{a+4}{2}}}Vol_g(B_{2\rho}).
\end{eqnarray}
By assumption, there exists an integer $s\geq0$ such that
$$Vol_g(B_{2\rho})\leq C\rho^s,$$
where $C$ is a constant independent of $\rho$.

Therefore we have
\begin{eqnarray}
&&\int_M\gamma^{a+4}|\vec{H}|^a|\nabla^\bot\vec{H}|^2d\mu_g  \nonumber
+(2m-1)\int_M\gamma^{a+4}|\vec{H}|^{a+4}d\mu_g \leq C(a)\rho^{s-\frac{a+4}{2}}.
\end{eqnarray}
Chosing $a$ such that $a>\max\{0, 2s-4\}$ and letting $\rho\to +\infty$ in the above inequality, we get
$$\int_M|\vec{H}|^a|\nabla^\bot\vec{H}|^2d\mu_g +(2m-1)\int_M|\vec{H}|^{a+4}d\mu_g =0,$$
which implies $\vec{H}=0$.

This completes the proof of theorem \ref{thm0}. \endproof
\subsection{Proof of theorem \ref{thm4}.}
From equation $(\ref{equ4})$ and the assumption that the sectional curvature of $(N, h)$ is small that $-\epsilon$ we see that
\begin{eqnarray*}
\Delta|\vec{H}|^2&=&2|\nabla^\bot\vec{H}|^2+2\langle\vec{H}, \Delta^\bot\vec{H}\rangle
\\&=&2|\nabla^\bot\vec{H}|^2+2\sum_{i=1}^m\langle B(A_{\vec{H}}e_i, e_i), \vec{H}\rangle-2\sum_{i=1}^m\langle R^N( e_i, \vec{H})e_i, \vec{H}\rangle
\\&\geq&2|\nabla^\bot\vec{H}|^2+2m|\vec{H}|^4+2m\epsilon|\vec{H}|^2
\\&\geq&2|\nabla^\bot\vec{H}|^2+2m\epsilon|\vec{H}|^2.
\end{eqnarray*}

Let $\gamma: M\to \R^+$ be a cut off function such that
$$\gamma=1\s on \s B_{\rho}, \gamma=0\s on\s M\setminus B_{2\rho},
\s and \s |\nabla\gamma|\leq\frac{C}{\rho},$$
for some constant $C$ independent of $\rho$. Here $B_\rho$ is a geodesic ball of radius $\rho$ on $M$. Then
\begin{eqnarray}\label{111}
&&-\int_M\nabla|\vec{H}|^2\nabla|\vec{H}|^a\gamma^2d\mu_g  \nonumber
\\&=&\int_M\Delta|\vec{H}|^2|\vec{H}|^a\gamma^2d\mu_g \nonumber
\\&\geq&2\int_M|\nabla^\bot\vec{H}|^2|\vec{H}|^a\gamma^2d\mu_g +2m\epsilon\int_M|\vec{H}|^{a+2}\gamma^2d\mu_g ,
\end{eqnarray}
where $a$ is a nonnegative constant.

On the other hand we have
\begin{eqnarray}\label{112}
&&-\int_M\nabla|\vec{H}|^2\nabla|\vec{H}|^a\gamma^2d\mu_g  \nonumber
\\&=&-\int_M2\langle\nabla^\bot\vec{H}, \vec{H}\rangle(a|\vec{H}|^{a-2}\langle\nabla^\bot\vec{H}, \vec{H}\rangle+2|\vec{H}|^a\gamma\nabla\gamma)d\mu_g  \nonumber
\\&\leq&-4\int_M\langle\nabla^\bot\vec{H}, \vec{H}\rangle|\vec{H}|^a\gamma\nabla\gamma d\mu_g  \nonumber
\\&\leq&2\int_M|\nabla^\bot\vec{H}|^2|\vec{H}|^a\gamma^2d\mu_g
+\frac{C}{\rho^2}\int_{B_{2\rho}\setminus B_\rho}|\vec{H}|^{a+2}d\mu_g.
\end{eqnarray}
From (\ref{111})-(\ref{112}), we get
\begin{eqnarray*}
2m\epsilon\int_{B_\rho}|\vec{H}|^{a+2}d\mu_g \leq\frac{C}{\rho^2}\int_{B_{2\rho}\setminus B_\rho}|\vec{H}|^{a+2}d\mu_g.
\end{eqnarray*}
Set $f(\rho)=\int_{B_\rho}|\vec{H}|^{a+2}d\mu_g$, from the above inequality we see that when $\rho$ is big enough, we have
$$f(\rho)\leq\frac{C}{\rho^2}f(2\rho).$$
This implies that $f(\rho)\leq\frac{C}{\rho^{2n}}f(2^n\rho)$, where $C$ is a constant independent of $\rho$. By assumption we have $f(2^n\rho)\leq C(1+2^{ns}\rho^s)$ for some positive integer $s$, as $\rho$ is big enough, hence $f(\rho)\leq\frac{C(1+2^{ns}\rho^s)}{\rho^{2n}}$. Let $2n>s$ one gets $\lim_{\rho\to\infty}f(\rho)=0,$ i.e. $\int_M|\vec{H}|^{a+2}d\mu_g=0$. Therefore, $\vec{H}=0$.

Furthermore, if $M$ is a hypersurface and the Ricci curvature of $(N, h)$ is smaller that $-\epsilon$, then
\begin{eqnarray*}
\Delta|\vec{H}|^2&=&2|\nabla^\bot\vec{H}|^2+2\langle\vec{H}, \Delta^\bot\vec{H}\rangle
\\&=&2|\nabla^\bot\vec{H}|^2+2\sum_{i=1}^m\langle B(A_{\vec{H}}e_i, e_i), \vec{H}\rangle-2\sum_{i=1}^m\langle R^N( e_i, \vec{H})e_i, \vec{H}\rangle
\\&=&2|\nabla^\bot\vec{H}|^2+2\sum_{i=1}^m\langle B(A_{\vec{H}}e_i, e_i), \vec{H}\rangle-2Ric(\vec{H}, \vec{H})
\\&\geq&2|\nabla^\bot\vec{H}|^2+2m|\vec{H}|^4+2\epsilon|\vec{H}|^2
\\&\geq&2|\nabla^\bot\vec{H}|^2+2\epsilon|\vec{H}|^2.
\end{eqnarray*}
Then by the same argument as before we can prove that $\vec{H}=0$.

This finishes the proof of theorem \ref{thm4}.
\endproof
\subsection{Proof of theorem \ref{thm2}.}
\proof We prove theorem \ref{thm2} in two steps.

 \textbf{Step 1.} At any point $x\in M$, $\vec{H}=0$ or $\nabla\vec{H}=0$.

 Let $\tilde{\gamma}: N\to \R^+$ be a cut off function defined as
$$\tilde{\gamma}=1\s on\s B_\rho,\tilde{\gamma}=0\s on\s N\setminus B_{2\rho}, \s |D\tilde{\gamma}|\leq\frac{C}{\rho}$$
and define $\gamma: M\to \R^+$ by $$\gamma:=\tilde{\gamma}\circ f.$$

By integral by parts we have
\begin{eqnarray}\label{ine1}
\int_M|\vec{H}|^a\langle\vec{H}, \Delta \vec{H}\rangle\gamma^2d\mu_g=&-&\int_M\vec{H}|^a|\nabla\vec{H}|^2\gamma^2d\mu_g
-2\int_M\vec{H}|^a\gamma\nabla\gamma\langle\nabla\vec{H}, \vec{H}\rangle d\mu_g \nonumber
\\&-&a\int_M|\vec{H}|^{a-2}\langle\nabla \vec{H}, \vec{H}\rangle^2\gamma^2 d\mu_g.
\end{eqnarray}
On the other hand, by assumption we have
\begin{eqnarray}\label{ine2}
\int_M|\vec{H}|^a\langle\vec{H}, \Delta \vec{H}\rangle\gamma^2d\mu_g\geq(\varepsilon-1)\int_M\vec{H}|^a|\nabla\vec{H}|^2\gamma^2d\mu_g.
\end{eqnarray}
From (\ref{ine1})-(\ref{ine2}), we have
\begin{eqnarray*}
&&\varepsilon\int_M\vec{H}|^a|\nabla\vec{H}|^2\gamma^2d\mu_g
\\&=&-2\int_M\vec{H}|^a\gamma\nabla\gamma\langle\nabla\vec{H}, \vec{H}\rangle d\mu_g \nonumber
-a\int_M|\vec{H}|^{a-2}\langle\nabla \vec{H}, \vec{H}\rangle^2\gamma^2 d\mu_g
\\&\leq&-2\int_M\vec{H}|^a\gamma\nabla\gamma\langle\nabla\vec{H}, \vec{H}\rangle d\mu_g,
\end{eqnarray*}
which implies by using Young's inequality that
\begin{eqnarray}
\varepsilon\int_M\vec{H}|^a|\nabla\vec{H}|^2\gamma^2d\mu_g\leq\frac{\varepsilon}{2}
\int_M\vec{H}|^a|\nabla\vec{H}|^2\gamma^2d\mu_g
+\frac{2}{\varepsilon}\int_M|\nabla\gamma|^2|\vec{H}|^{a+2}d\mu_g.
\end{eqnarray}
Therefore
\begin{eqnarray*}
\int_{f^{-1}(B_{\rho})}|\vec{H}|^a|\nabla\vec{H}|^2d\mu_g
&\leq&\int_M\vec{H}|^a|\nabla\vec{H}|^2\gamma^2d\mu_g
\\&\leq&\frac{4}{\varepsilon^2}
\frac{C}{\rho^2}\int_{f^{-1}(B_{2\rho})}|\vec{H}|^{a+2}d\mu_g.
\end{eqnarray*}
Let $\rho\to \infty$ in this inequality we get
\begin{eqnarray}
\int_M|\vec{H}|^a|\nabla\vec{H}|^2d\mu_g\leq 0,
\end{eqnarray}
and hence $\vec{H}=0$ or $\nabla\vec{H}=0$.

\textbf{Step 2.} $\nabla\vec{H}=0$ implies $\vec{H}=0$.

Now let $x\in M$ such that $\nabla\vec{H}(x)=0$.
We choose an orthonormal basis $\{e_i, i=1,...,m\}$ of $T_xM$ and
an orthonormal basis $\{\nu^\alpha, \alpha=1,...,t\}$ of $(T_xM)^\bot$.

Computing directly one gets
\begin{eqnarray}
0=\langle\nabla_{e_i} \vec{H}, e_j\rangle=\sum_\alpha H_\alpha B_\alpha(e_i, e_j),
\end{eqnarray}
for any $1\leq i,j\leq m$. Taking trace over this equality we get
$$\sum_\alpha (H_\alpha)^2=|\vec{H}|^2=0.$$
Therefore $\vec{H}=0$.

Combing step 1 and 2, we see that if $f: (M, g)\to (N^{m+t}, h)$ is a biharmonic submanifold satisfying
the assumptions of theorem \ref{thm2}, then $\vec{H}=0$.
\endproof
\subsection{Proof of theorem \ref{thm3}.}
The proof is quite similar as before. For the convenience of readers we give the details.

 \proof
 Let $\gamma: M\to \R^+$ be a cut off function such that
$$\gamma=1\s on \s B_{\rho}, \gamma=0\s on\s M\setminus B_{2\rho},
\s and \s |\nabla\gamma|\leq\frac{C}{\rho},$$
for some constant $C$ independent of $\rho$. Here $B_\rho$ is a geodesic ball of radius $\rho$ on $M$.

By integral by parts we get
\begin{eqnarray}
\int_M|\vec{H}|^a\langle\vec{H}, \Delta \vec{H}\rangle\gamma^2d\mu_g=&-&\int_M\vec{H}|^a|\nabla\vec{H}|^2\gamma^2d\mu_g
-2\int_M\vec{H}|^a\gamma\nabla\gamma\langle\nabla\vec{H}, \vec{H}\rangle d\mu_g \nonumber
\\&-&a\int_M|\vec{H}|^{a-2}\langle\nabla \vec{H}, \vec{H}\rangle^2\gamma^2 d\mu_g.
\end{eqnarray}
On the other hand, by assumption we have
\begin{eqnarray}
\int_M|\vec{H}|^a\langle\vec{H}, \Delta \vec{H}\rangle\gamma^2d\mu_g\geq(\varepsilon-1)\int_M\vec{H}|^a|\nabla\vec{H}|^2\gamma^2d\mu_g.
\end{eqnarray}
Hence we have, by combining the above two inequalities, that
\begin{eqnarray*}
&&\varepsilon\int_M\vec{H}|^a|\nabla\vec{H}|^2\gamma^2d\mu_g
\\&=&-2\int_M\vec{H}|^a\gamma\nabla\gamma\langle\nabla\vec{H}, \vec{H}\rangle d\mu_g \nonumber
-a\int_M|\vec{H}|^{a-2}\langle\nabla \vec{H}, \vec{H}\rangle^2\gamma^2 d\mu_g
\\&\leq&-2\int_M\vec{H}|^a\gamma\nabla\gamma\langle\nabla\vec{H}, \vec{H}\rangle d\mu_g,
\end{eqnarray*}
which implies by using Young's inequality that
\begin{eqnarray}
\varepsilon\int_M\vec{H}|^a|\nabla\vec{H}|^2\gamma^2d\mu_g\leq\frac{\varepsilon}{2}
\int_M\vec{H}|^a|\nabla\vec{H}|^2\gamma^2d\mu_g
+\frac{2}{\varepsilon}\int_M|\nabla\gamma|^2|\vec{H}|^{a+2}d\mu_g.
\end{eqnarray}
Noth that
$$\int_{B_{2\rho}}|\vec{H}|^{a+2}d\mu_g\leq\int_M|\vec{H}|^{a+2}d\mu_g,$$
therefore one gets
\begin{eqnarray}
\int_{B_{\rho}}|\vec{H}|^a|\nabla\vec{H}|^2d\mu_g\leq\frac{4}{\varepsilon^2}
\frac{C}{\rho^2}\int_M|\vec{H}|^{a+2}d\mu_g.
\end{eqnarray}
Let $\rho\to \infty$ in this inequality we get
\begin{eqnarray}
\int_M|\vec{H}|^a|\nabla\vec{H}|^2d\mu_g\leq 0,
\end{eqnarray}
and hence $\vec{H}=0$ or $\nabla\vec{H}=0$.

Now using the same procedure as step 2 of the proof of theorem \ref{thm2}, we can prove that $\vec{H}=0$.
\endproof

  \textbf{Acknowledgement.} The author would like to express his appreciation to Professor Guofang Wang for his suggestions on this paper. He is also very appreciated with Dr. Shun Maeta for pointing out a mistake in the manuscript.
 {}
\vspace{1cm}\sc
Yong Luo

School of mathematics and statistics,

Wuhan university, Hubei 430072, China

{\tt yluo@amss.ac.cn}

\vspace{1cm}\sc

\end{document}